\documentclass[11pt]{article}
\usepackage{amsmath}
\usepackage{amsthm}
\usepackage{amssymb}

\newtheorem{theo}{Theorem}[section]

\newtheorem{prop}[theo]{Proposition}
\newtheorem{lemma}[theo]{Lemma}

\newcommand{\CC}{{\cal C}}

\begin{document}
\date{}

\title{
Explicit expanders of every degree and size
}

\author{Noga Alon
\thanks
{Department of Mathematics, Princeton University,
Princeton, NJ 08544, USA and
Schools of Mathematics and
Computer Science, Tel Aviv University, Tel Aviv 69978,
Israel.  
Email: {\tt nogaa@tau.ac.il}.  
Research supported in part by
NSF grant DMS-1855464, ISF grant 281/17, BSF grant 2018267
and the Simons Foundation.}
}

\maketitle
\begin{abstract}
An $(n,d,\lambda)$-graph is a $d$ regular graph on $n$ vertices in which
the absolute value of any nontrivial eigenvalue is at most 
$\lambda$.
For any constant $d \geq 3$, $\epsilon>0$ and all sufficiently 
large $n$ we show that there is a deterministic poly(n) time algorithm
that outputs an $(n,d, \lambda)$-graph (on exactly $n$ vertices)
with $\lambda \leq 2 \sqrt{d-1}+\epsilon$.
For any $d=p+2$ with $p \equiv 1 \bmod 4$
prime and all sufficiently  large $n$, we describe a strongly 
explicit construction of an $(n,d, \lambda)$-graph (on exactly
$n$ vertices) with $\lambda \leq \sqrt {2(d-1)} + \sqrt{d-2} +o(1)
(< (1+\sqrt 2) \sqrt {d-1}+o(1))$, with the $o(1)$ term tending to
$0$ as $n$ tends to infinity.  For every $\epsilon >0$, 
$d>d_0(\epsilon)$ and $n>n_0(d,\epsilon)$ we  present a strongly explicit
construction of an $(m,d,\lambda)$-graph with $\lambda < (2+\epsilon) 
\sqrt d$ and $m=n+o(n)$. All constructions are obtained by starting with
known ones of Ramanujan or nearly Ramanujan graphs, modifying
or packing them in an appropriate way. The spectral analysis relies on
the delocalization of eigenvectors of regular 
graphs in cycle-free neighborhoods.
\end{abstract}

\section{Introduction}

An $(n,d,\lambda)$-graph is a $d$-regular graph on $n$ vertices 
in which the absolute value of every nontrivial eigenvalue is at most
$\lambda$. This notation was introduced by the author in the early
90s motivated by the fact that such graphs in which $\lambda$ is much
smaller than $d$ exhibit strong expansion and quasi-random properties,
see \cite{Al1}, \cite{AC}, \cite{KS}.

It is well known (see \cite{Al1}, \cite{Ni}, \cite{Fr1})
that if an $(n,d,\lambda)$-graph exists then $\lambda \geq 2\sqrt {d-1}
-O(1/\log^2 n)$. An $(n,d,\lambda)$-graph is called (two-sided) 
{\em Ramanujan} if
$\lambda \leq 2\sqrt {d-1}$.

Lubotzky, Phillips and Sarnak \cite{LPS}, and 
independently Margulis \cite{Ma} proved that
for every prime $p$ which is $1$ modulo $4$ there are infinite families
of $d$-regular Ramanujan graphs. Friedman \cite{Fr2} (see also
\cite{Bo} for a simpler proof) proved the existence of near
Ramanujan graphs of every degree and every (large) admissible size. Indeed,
establishing a conjecture of the present author he proved that a random
$d$-regular graph on $n$ vertices is, with high probability, an 
$(n,d,\lambda)$-graph for $\lambda=2\sqrt{d-1}+o(1)$, where the $o(1)$-term
tends to zero as $n$ tends to infinity.

For applications, however, (see, e.g., \cite{HLW} and its references for 
many of those) it is desirable to have explicit constructions
of such graphs.  It is also sometime desirable to have explicit
constructions with specified degrees and number of vertices,
(see, for example, \cite{MRSV} for a recent example). 
A construction is called {\em explicit} if there is a 
deterministic polynomial time algorithm that, given $n$ and $d$,
produces an $(n,d,\lambda)$-graph   
(or an $(n(1+o(1)), d, \lambda)$-graph). It is {\em strongly explicit }
if the adjacency list of any given vertex can be produced in time
$polylog(n)$.  The construction of \cite{LPS}, and that of \cite{Ma}
are strongly explicit
\footnote {Though they require to find a large prime in a
prescribed range. This can be done efficiently using randomization,
but can also be avoided. More details appear in Section 2}, 
providing Cayley graphs of $SL(2,F_q)$, but work
only for degrees that are $p+1$ for primes $p \equiv 1 \bmod 4$ and for
numbers of vertices that are of the form $q(q^2-1)/2$ for  primes
$q$ which are $1$ modulo $4$ so that $p$ is a quadratic residue modulo
$p$. Morgenstern \cite{Mo} gave a strongly explicit construction
for every degree which is a prime power plus 1, but the possible numbers
of vertices obtained  are sparser.  An observation in \cite{CM} provides
strongly  explicit families of $(n,d,\lambda)$-graphs
with $\lambda \leq O(d^{0.525})$ for infinitely many values of $n$
(but not for every $n$). Similarly, the method in 
\cite{RVW} and its improvement in
\cite{BaTs} provide strongly explicit families with $\lambda \leq
O(d^{1/2+o(1)})$ (for infinitely many, but not for all $n$). The results
of \cite{MSS} together with those of \cite{Co} and an observation of
Srivastava (cf. \cite{MOP}) give explicit, but not strongly explicit
$(n,d,\lambda)$-graphs for all admissible $d$ and $n$ with
$\lambda \leq 4 \sqrt {d-1}$. In a recent work of Mohanty, O'Donnell
and Paredes \cite{MOP} 
the authors describe an explicit (but not strongly explicit)
construction of $(n,d,\lambda)$-graphs for every $d$, where
$\lambda = 2\sqrt {d-1}+o(1)$ and the $o(1)$-term tends to $0$ as
$n$ tends to infinity. This, again, works for infinitely many values
of $n$, but not for all $n$.

In the present short paper 
we describe improved explicit and strongly explicit
constructions of near Ramanujan graphs of all degrees and (large) number of
vertices. The first result is a
(slightly improved version of an) 
observation I mentioned in several lectures in the 90s that, as far
as I know, has never appeared in print. Although it is very simple,
the parameters it provides are far better than the ones obtained from
the constructions in \cite{CM}, \cite{RVW}, \cite{BaTs}, and I therefore
decided to include it here.
\begin{prop}
\label{p11}
For every degree $d$ there is a strongly explicit constructions of
$(n,d,\lambda)$-graphs where $\lambda \leq (2+o_d(1))\sqrt {d}$,
the $o_d(1)$-term tends to zero as $d$ tends to infinity,
and the possible values of $n$ form a sequence in which the ratio
between consecutive terms tends to $1$.
\end{prop}
Note that this means that for every desired number of vertices $n$ and any
desired degree $d$, there is a strongly explicit construction of
an $(n(1+o_n(1)),d, \lambda)$-graph with $\lambda \leq (2+o_d(1)) 
\sqrt d$. Here the term $o_n(1)$ tends to zero as $n$ tends to infinity
and the $o_d(1)$-term tends to zero as $d$ tends to infinity.

The next result provides strongly explicit constructions of 
$(n,d,\lambda)$ graphs for degrees $d=p+2$ with $p$ being a prime
congruent to $1$ modulo $4$, for any desired (large) number of vertices.
\begin{theo}
\label{t12}
For any prime $p \equiv 1 \bmod 4$ and every sufficiently large $n$
there is a strongly explicit construction of an $(n,d,\lambda)$-graph
(on exactly $n$ vertices), where $d=p+2$ and 
$\lambda \leq \sqrt{2(d-1)}+\sqrt{d-1} +o(1)
< (1+\sqrt 2) \sqrt {d-1} +o(1)$, and the $o(1)$-term tends to
zero as $n$ tends to infinity.
\end{theo}
It is worth noting that here we allow to have at most one loop in
every vertex, with the convention that a loop adds one to the degree
(otherwise we must have an even number of vertices as the degree of
regularity is odd). For even $n$ we can replace the loops by a
matching with no loss in the spectral estimate.

If an explicit, rather than strongly  explicit construction suffices,
we can combine a variant of 
our method with the new result of \cite{MOP} to get
the following.

\begin{theo}
\label{t13}
For every degree $d$, every $\epsilon$ and all sufficiently 
large $n \geq n_0(d,\epsilon)$, where $nd$ is even,
there is an explicit
construction of an $(n,d,\lambda)$-graph with $\lambda \leq 2\sqrt{d-1}
+\epsilon$.
\end{theo}

The construction in the proof of  Proposition \ref{p11} is a 
simple packing of known Ramanujan graphs on the same set of vertices. 
A crucial point is that these constructions are all Cayley graphs of the
same group, so one can simply take a union of the 
corresponding generating sets.
The proofs of Theorem \ref{t12} and \ref{t13} require more work.
Here too the idea is to start from a known  Ramanujan or nearly Ramanujan
graph and modify it in an appropriate way. In the proof of
Theorem \ref{t12} we add
vertices connected to arbitrary disjoint sets of neighbors, adding loops
(or a matching)
to keep the graph regular. The eigenvalues are then estimated by their
variational definition.
In the construction for Theorem \ref{t13}
we omit carefully chosen vertices from a given near-Ramanujan graph
and add a matching between their neighbors to maintain regularity.
A crucial point in the spectral analysis here is the delocalization of
the eigenvectors of the graphs obtained, which is based on the  absence
of short cycles in the neighborhoods of the omitted vertices.

The rest of this paper is organized as follows. In Section 2 we describe the
strongly explicit constructions, including the proofs of Proposition
\ref{p11} and Theorem \ref{t12}. In Section 3 we present the proof of
Theorem \ref{t13}. The final Section 4 contains some concluding remarks
and open problems.

\section{Strongly explicit constructions}

The basic construction we describe here requires the ability to
find efficiently a large prime in a prescribed range. It is well
known that this can be done efficiently by a randomized algorithm,
and can also be done deterministically assuming some standard
(open) number-theoretic conjectures 
about the gap between consecutive primes. Since
this is the only non-deterministic part of the construction, we
call it a $p$-strongly explicit construction (where $p$ stands for 
prime). This construction is described in the first subsection.
We then show how it can be replaced by a totally strongly
explicit construction. To do so, we first include a subsection
presenting the (known)  description of the construction of
\cite{LPS}, \cite{Ma} as Cayley graphs of Quaternions over $Z_m$.
We proceed with a proof of Theorem \ref{t12} with a $p$-strongly
explicit construction, followed by its modification to a strongly
explicit one. 

\subsection{The basic construction}
\label{basic}

We start with the simple proof of 
Proposition \ref{p11}, with a $p$-strongly explicit construction. 
It is based on the
fact that if $G_i=(V,E_i)$, $i \in I$, are graphs on the same set of
vertices $V$, where $G_i$ is an $(n,d_i,\lambda_i)$-graph, then their
union $G=(V, \cup_i E_i)$ (considered as a multigraph in case the sets
$E_i$ are not pairwise disjoint), is an $(n, \sum_i d_i, \sum_i \lambda_i)$
graph. This is a simple consequence of the variational definition of the
eigenvalues. The Ramanujan graphs in \cite{LPS} or \cite{Ma} are Cayley
graphs of the group $SL(2,F_q)$ of the two by two matrices with determinant
$1$ over the finite field $F_q$, modulo 
its normal subgroup consisting of the
identity $I$ and $-I$.  The degree can be $1$ plus 
any prime $p$ congruent to $1$
modulo $4$, where $q$ is also a prime congruent to $1$ modulo $4$, and
$p$ is a quadratic residue modulo $q$. Note that by quadratic reciprocity
this is equivalent to $q$ being a quadratic residue modulo $p$.

Given a desired degree $d=d_1$, 
let $p_1 $ be the largest prime congruent to
$1$ modulo $4$ and satisfying
$p_1+1 \leq d_1$. Put $d_2=d_1-p_1-1$. If $d_2>4$ 
let $p_2$ be the largest prime
congruent to $1$ modulo $4$ which satisfies $p_2+1 \leq d_2$ and put
$d_3=d_2-p_2-1$. Continuing in this manner we get primes $p_1, \ldots ,
p_s$ as above so that $(p_1+1)+(p_2+1)+ \cdots +(p_s+1) \leq d$
where $y=d-((p_1+1)+(p_2+1)+ \cdots +(p_s+1)) \leq 4$. Let $q$
be a prime  congruent to $1$ modulo $4$ which is a quadratic residue 
modulo each $p_i$ (for example, any $q$ which is $1$  modulo each $p_i$
will do). Let $V$ be the
set of elements of
$SL(2,F_q)$. For each $i$ let $G_i$ be the $(p_i+1)$-regular 
Ramanujan 
Cayley graph of $SL(2,F_q)$ described in \cite{LPS}, and let
$X_i$ be its (symmetric) set of generators. Let $G'$ be the 
Cayley graph of $SL(2,F_q)$ whose set of generators consists of the union
of all sets $X_i$. Then $G'$ is $(d-y)$-regular, where $0 \leq y \leq 4$.
If $y=0$ let $G$ be $G'$. If $y=1$ add to the set of generators the 
matrix $M$ with rows $(0,1)$ and $(-1,0)$ (which is of order $2$). If
$y=2$ add an arbitrary generator and its inverse, if $y=3$ add such
a generator, its inverse and $M$, and if $y=4$ add an arbitrary set of
two generators and their inverses. In each of these cases the resulting
graph $G$ is a $d$-regular Cayley graph of $SL(2,F_q)$. By the known
results about the distribution of primes in arithmetic progressions
each prime $p_i$ is much smaller than $p_{i-1}$ as long as
$p_{i-1}$ is large. In fact, by \cite{BHP} it follows
that $p_{i}=O(p_{i-1}^{0.525})$.
Therefore, the resulting
graph $G$ is an $(n,d,\lambda)$-graph for $n=q(q^2-1/2$ with
$\lambda\leq (2+o_d(1))\sqrt d$, where the $o(1)$-term tends to
zero as $d$ tends to infinity. Note that it is not difficult to ensure,
if so desired, 
that the graph $G$ is simple: we just have to ensure the chosen
primes are distinct. This is automatically the case whenever $d_i$ is
still large, and if needed we can stop when $d_i$ becomes small and add
arbitrary additional generators and their inverses, together with 
$M$ if $d$ is odd.  Alternatively, if we have
to repeat the same prime several times, we can take the corresponding
generating set for this prime 
and conjugate it to get an isomorphic graph with 
different generators.
The known results about the
distribution of primes in arithmetic progressions imply also that 
for each choice of the primes $p_i$ the possible choices for the prime
$q$ suffice to ensure that the sequence of possible values for the number
of vertices $n$ of the graph is one in which the ratio between consecutive
terms tends to $1$ as $n$ tends to infinity. This completes the proof
of the proposition (with a $p$-strongly explicit construction
resulting from the need to find the required large prime $q$).  
\hfill $\Box$
\vspace{0.2cm}

\subsection{Ramanujan graphs as Cayley graphs of quaternions}
\label{quaternions}

In this subsection we present the known description of the LPS Ramanujan
graphs as Cayley graphs of quaternions. The proof these are
Ramanujan graphs  appears
(somewhat implicitly)
in \cite{Lu}. 

Let $p$ be a prime congruent to $1$ modulo $4$, and
let $A=A(p)$ be the set of all integral solutions $(a_0,a_1,a_2,a_3)$ 
of the equation $a_0^2+a_1^2+a_2^2+a_3^2=p$ 
where $a_0$ is positive odd,
and all other $a_i$ are even. By a well known  result of Jacobi 
there are exactly $p+1$  such vectors.
Let $m$ be odd, relatively prime to $p$, and assume further that $p$ is
a square in $Z_m^*$. 
let $Q(m)$ be the factor group of the 
multiplicative group of the quaternions 
over $Z_m$ whose norm
is a square in $Z_m^*$, modulo its normal subgroup consisting of
the scalars $Z_m^*$. Thus the elements of $Q(m)$ are all 
quaternions $x_0+x_1 i+x_2 j+x_3 k$ where 
$x_0^2+x_1^2+x_2^2+x_3^2 \in (Z_m^*)^2$ and two such elements are
identified if one is a multiple of the other by a scalar.
Finally, let $H=H(p,m)$ be the Cayley graph of $Q(m)$ with 
the generating set 
$$\{a_0+a_1i+a_2j+a_3k: (a_0,a_1,a_2,a_3) \in A(p)\}.$$ 
The following result is proved (somewhat implicitly) in \cite{Lu},
see pages 95-97.
\begin{theo}[\cite{Lu}]
\label{t91}
For every $p$ and $m$ as above $H=H(p,m)$ is a non-bipartite
$(p+1)$-regular Ramanujan  graph, that is, the absolute value of each
of its eigenvalues besides the top one is at most $2 \sqrt p$.
\end{theo}

\subsection{The proof of Proposition \ref{p11}}

In the construction here we will 
start with the graphs $Q(p,m)$ with $p \equiv
1 \bmod 4$ a prime and $m=q_1^s q_2^t$, where $s,t \geq 1$ and
$q_1,q_2 $ are distinct primes, each being $1 \bmod 4p$. For each
fixed $p$ as above, the known results about the Linnik problem 
(see \cite{HB}) 
imply that there are $q_1,q_2$ as above, each being at most a polynomial
in $p$.  It is not difficult to check, using Hensel's Lemma and the 
Chinese Remainder Theorem, that the number of vertices
of $H(p,q_1^s q_2^t)$ is 
$$
Q(q_1,q_2,s,t)=q_1^{3(s-1)} q_2^{3(t-1)} \frac{q_1(q_1-1)(q_1+1)}{2} 
\frac{q_2(q_2-1)(q_2+1)}{2}.
$$
Indeed, by Hensel's Lemma, for elements $x_0,x_1,x_2,x_3$ of $Z_m$ 
the norm $x_0^2+x_1^2+x_2^2+x_3^2$ is a square in $Z^*_m$ if and
only if it is a square in $Z^*_{q_1}$ and in $Z^*_{q_2}$. Since
each $q_i$ is $1 \bmod 4$, $-1$ is a quadratic residue implying
that the number of solutions of $y_1^2+y_2^2=0$ in $Z_{q_i}$ is
$2q_i-1$. For each nonzero $b$ in $Z_{q_i}$
the number of solutions of $y_1^2+y_2^2=b$ (in $Z_{q_i}$) is the same
as the number of solutions  of $y^2-z^2(=(y-z)(y+z))=b$ , which is 
$q_i-1$. This shows that the number of solutions of
$x_0^2+x_1^2+x_2^2+x_3^2=b$ for any nonzero $b \in Z_{q_i}$ is
$$
2 (2q_i-1)(q_i-1)+(q_i-2)(q_i-1)^2=(q_i-1)q_i(q_i+1).
$$
(These include $(2q_i-1)(q_i-1)$ solutions with
$x_0^2+x_1^2=0$ and $x_2^2+x_3^2=b$, $(2q_i-1)(q_i-1)$ ones
with $x_0^2+x_1^2=b$ and $x_2^2+x_3^2=0$, and $(q_i-1)^2$ solutions
for each of the $q_i-2$ possibilities $x_0^2+x_1^2=b_1$ and 
$x_2^2+x_3^2=b_2$ with $b_1+b_2=b$ and $b_1,b_2 \not \in \{0,b\}$.)
Therefore, the number of elements over $Z_{q_i}$ whose norm is
a nonzero square in $Z_{q_i}$ is 
$$
\frac{q_i-1}{2} (q_i-1)q_i(q_i+1)=\frac{q_i(q_i-1)^2(q_i+1)}{2}.
$$
By the Chinese remainder Theorem there are 
$$
\frac{q_1(q_1-1)^2(q_1+1)}{2} 
\frac{q_2(q_2-1)^2(q_2+1)}{2}.
$$
elements $(x_0,x_1,x_2,x_3)$ in $Z_{q_1q_2}$ so that
$x_0^2+x_1^2+x_2^2+x_3^2$ is a square in $Z^*_{q_1q_2}$, and
by Hensel's Lemma  each of them provides $q_1^{4(s-1)}q_2^{4(t-1)}$
quaternions over $Z_{q_1^sq_2^t}$ whose norm is a square in
$Z^*_{q_1^sq_2^t}$. To get the number of vertices of the graph we
just have to divide by the cardinality of $Z^*_{q_1^sq_2^t}$ which
is $q_1^{s-1}q_2^{t-1}(q_1-1)(q_2-1)$, obtaining the required 
number of vertices.  Note also that by this description it is
easy to number the vertices of the graph. (For our application here it
is in fact enough to number a constant fraction of them. For fixed
$q_1,q_2$ this can be done, for example, by numbering all vectors
$(1,x_1,x_2,x_3) \in Z_{q_1^sq_2^t}$ with $x_1,x_2,x_3$ divisible by
$q_1q_2$, lexicographically).

We next show that for every fixed distinct 
primes $q_1,q_2$, the ratio between
consecutive elements in the set of integers
$\{Q(q_1,q_2,s,t): s,t \geq 1 \}$ tends to $1$ as the elements grow.
\begin{lemma}
\label{l32}
Let $q_1,q_2$ be distinct primes. Then for every large integer
$n$ there are positive integers $s,t$ so that
$n \leq Q(q_1,q_2,s,t) \leq n+o(n)$.
\end{lemma}

\noindent
{\bf Proof:\,} The constant 
$\alpha=\frac{\log {q_1}}{\log {q_2}}$ is irrational. Therefore, by
the equidistribtion theorem (in fact, by a special case that follows
easily from the pigeonhole principle), for every $\delta>0$ there
is an integer $k_1=k_1(\alpha)$ so that 
$0<k_1 \alpha \bmod 1 < \delta$.
It follows that for every $\mu>0$ there are  
integers $k_1,k_2$ such that
$$
1 \leq \frac{q_1^{k_1}}{q_2^{k_2}} \leq q_2^{\delta} \leq (1+\mu).
$$
This implies that for every $s,t \geq \max\{k_1,k_2\}$ the ratio
between $Q(q_1,q_2,s,t)$ and $Q(q_1,q_2,s-k_1,t-k_2)$ is between
$1$ and $(1+\mu)^3$, implying the desired result. \hfill $\Box$
\vspace{0.2cm}

\noindent
The proof of Proposition \ref{p11} now proceeds exactly as in
subsection \ref{basic}, using the description of the LPS graphs 
serving as the building blocks as given in subsection
\ref{quaternions}. Since here $q_1,q_2$ are constants, there is no
need to find any large primes for the construction, providing a
strongly explicit construction for every fixed degree.

\subsection{The proof of Theorem \ref{t12}}

We first describe a $p$-strongly-explicit construction, starting,
again, with
the graphs of \cite{LPS}. Recall that the vertex sets of these graphs
is the set of matrices in $SL(2,F_q)$ where each matrix $A$ is identified
with $-A$. It is easy to number the vertices starting with the matrices
$(a_{ij})$ with $a_{11} \neq 0$ and ordering 
them according to the lexicographic
order of the elements $(a_{11},a_{12},a_{13})$ where $a_{14}$ is chosen to
ensure that the determinant is $1$ 
(which is always possible as $a_{11} \neq 0$). Here $1 \leq a_{11} \leq
(q-1)/2$, as we identify each matrix $A$ with $-A$. The first matrices
are the $q^2$ matrices with $a_{11}=1$, then those with
$a_{11}=2$, and so on. 
(The remaining $q(q-1)/2$ matrices with $a_{11}=0$ can appear 
last in our order according to the lexicographic order of 
$(a_{12},a_{24})$, but this will play no real role in the
construction.)

Given the desired number $n$ of vertices, and given the degree $d=p+2$
with $p$ as in the theorem, let $q$ be the largest prime which
is $1$ modulo $4$, is a quadratic residue modulo $p$ and satisfies
$|SL(2,F_q)|=m_q=q(q^2-1)/2 \leq n$. Put $m=m_q$, let $H$ be the Ramanujan 
$(p+1)=(d-1)$-regular graph of \cite{LPS} whose vertex set $V$ is the set
of elements of $SL(2,F_q)$ numbered as described above.
By the known results about the distribution
of primes in progressions $n-m=o(m)$. Put $r=n-m$ and let
$R$ be a set of $r$ additional vertices $u_1,u_2, \ldots, u_r$.
Connect each vertex $u_i$ to the vertices numbered
$(i-1)d+1,(i-1)d+2, \ldots ,id$ of $H$. Finally add
a loop to each remaining vertex of $H$ to make the graph regular.
This is the desired graph  $G$. It is clearly $d=p+2$-regular.
(If $n$ is even and we do not want loops we can replace them
by a matching between consecutive pairs of vertices,
saturating all non-neighbors of the $r$ new vertices). 

It is clear that the construction above is strongly explicit. 
To complete the proof
it remains to show that the absolute value of any nontrivial eigenvalue
of $G$ is at most $\sqrt{2(p+1)}+ \sqrt{p}+o(1)$.  
We proceed with a proof of this
fact. By the variational definition of the nontrivial eigenvalues of
$G$ this is equivalent to
showing that for every real function $f(u)$ on the set of vertices
$U=V \cup R$ of $G$ satisfying  $\|f \|_2^2=1$
and $\sum_{u \in U} f(u)=0$
\begin{equation}
\label{e31}
|f^t A_G f| \leq
\sqrt{2(p+1)}+ \sqrt{p}+o(1)
\end{equation}
where $A_G$ is the adjacency matrix of $G$. Let $W \subset V$ 
denote the set of all $(p+2)r$ neighbors of $R$, put $L=V-W$, and let
$E_R$ denote the set of all edges between $R$ and $W$. Thus $E_R$ is a 
collection of pairwise vertex disjoint stars, each having $(p+2)$ leaves.
The adjacency matrix of $G$ can be written as a sum
$A_G=A_H+A_R+A_L$, where $A_H$ is the adjacency matrix of the 
Ramanujan graph $H$ (with the additional isolated vertices
of $R$), $A_R$ is the adjacency matrix of the graph
$(U,E_R)$, and $A_L$ is the adjacency matrix of the graph on $U$ whose
edges are 
the loops on the vertices of $L$ (or the added matching on them, if 
we have chosen not to add loops).
Therefore
\begin{equation}
\label{e32}
f^t A_G f = f^t A_H f+f^t A_R f+f^t A_L f.
\end{equation}
We proceed to bound  each of these terms.

By Cauchy-Schwarz 
$$
|\sum_{u \in R} f(u)|^2 \leq |R| \sum_{u \in R} f^2(u) \leq |R|=o(m).
$$
Since $\sum_{u \in U} f(u)=0$, this implies that 
$|\sum_{u \in V} f(u)|^2= |\sum_{u \in R} f(u)|^2 =o(m)$

Let $g$ be the trivial normalized eigenvector of $H$, that is, the vector
given by $g(v)=1/\sqrt m$ for all $v \in V$. Expressing  the restriction
$f'$ of $f$ to $V$ as a linear combination of $g$ and a unit vector $h$
orthogonal to it, we get $f'=bg+ch$, where $\sum_{v \in V} h(v)=0$,
$b^2+c^2=1$ and $b^2=|\sum_{u \in V} f(u)|^2/m=o(1)$.
Since $H$ is a Ramanujan graph, $|h^t A_H h| \leq 2 \sqrt p.$
Therefore
\begin{equation}
\label{e33}
|f^t A_H f| = |(f')^t A_H f'| \leq b^2 (p+1) +c^2 2 \sqrt p
\leq 2 \sqrt p +o(1).
\end{equation}
Clearly
\begin{equation}
\label{e34}
|f^t A_L f| \leq \sum_{v \in L} f^2(v).
\end{equation}
Indeed this is an equality if there are loops and an inequality in
case  a matching has been added.

For bounding the absolute value of $f^t A_R f$ observe that for every 
positive $x$
$$
|f^t A_R f| =2 |\sum_{uv \in E_R} f(u) f(v) | 
$$
\begin{equation}
\label{e35}
\leq \sum_{u \in R, v \in W, uv \in E_R} (\frac{f^2(u)}{x} + x f^2(v))
=\frac{p+2}{x} \sum_{u \in R} f^2 (u) + x \sum_{v \in w} f^2 (v).
\end{equation}
Combining (\ref{e32}),(\ref{e33}),(\ref{e34}) and (\ref{e35}) we conclude
that for every positive real $x$
\begin{equation}
\label{e36}
|f^t A_G f| \leq (2\sqrt p+1) \sum_{v \in L} f^2 (v) +
(2 \sqrt p+x)\sum_{v \in W} f^2 (v) + \frac{p+2}{x} \sum_{v \in R}
f^2(v) + o(1).
\end{equation}
Choosing $x=\sqrt{2p+2}-\sqrt{p}$ (which is at least $1$) and substituting
in (\ref{e36}) we finally get
$$
|f^t A_G f| \leq (\sqrt{2p+2}+\sqrt p) \sum_{u \in U} f^2 (u)+o(1)
=(\sqrt{2p+2}+\sqrt p) +o(1).
$$
This establishes (\ref{e31}) and completes the proof (with a
$p$-strongly explicit construction). The conversion  to a strongly
explicit construction proceeds just as in the proof of Proposition
\ref{p11}, based on the results in subsection \ref{quaternions}. 
Note that as mentioned in that subsection the description there
provides a simple efficient  way to number enough vertices of each
graph $H(p,q_1^sq_2^t)$ and by Lemma \ref{l32} we can start by
finding efficiently appropriate $s,t$ using binary search. 
\hfill $\Box$

\section{Explicit constructions}

In this section we present the proof of Theorem \ref{t13}. We start
with some preliminary lemmas.
\begin{lemma}
\label{l31}
Let $G=(V,E)$ be a $d$-regular graph on $n$ vertices, where
$d \geq 3$, and suppose that
the $2r+4$-neighborhood of any vertex in it contains at most one cycle,
where $r \leq \log_{d-1} n$.
Then there is a subset $U \subset V$ of vertices of $G$ satisfying the 
following.
\begin{enumerate}
\item
$|U| \geq \frac{n}{2d^{2r+3}}$
\item
The $(r+1)$-neighborhood of any vertex in $U$ contains no cycle.
\item
The distance between any two vertices in $U$ is at least $2r+3$.
\end{enumerate}
Such a set $U$ can be found in polynomial time.
\end{lemma}
\vspace{0.1cm}

\noindent
{\bf Proof:}\,  Let $\CC$ denote the collection of all cycles of length
at most $2r+4$ in $G$. Note that the distance between any two members 
$C_1,C_2$ of
$\CC$ is larger than $2r+4$, since otherwise there is a vertex $v$
within distance at most $r+2$ of both cycles $C_i$, and then its
$2r+4$-neighborhood contains both cycles, contradiction. 
The $r+2$ neighborhood of each cycle $C \in \CC$ contains no
other cycle besides $C$, as it is contained in the
$2r+4$ neighborhood of any vertex on the cycle. Thus the number of edges 
spanned by each such neighborhood is at most the number of vertices
in it. As the neighborhoods are vertex disjoint, the total number
of edges in all these neighborhoods together is at most the number
of vertices of $G$ which is $n$. It follows that by omitting all
vertices in the $(r+1)$-neighborhoods of all members of $\CC$, at most
$n$ edges are omitted, and as $G$ has $nd/2$ edges and $d \geq 3$
at least $n/2 $ edges, and hence at least $n/2d$ vertices
have not been omitted. Let $Z$ be the set
of non-omitted vertices. Note that the $(r+1)$-neighborhood of any vertex
in $Z$ contains no cycle (as if it contains a cycle, it contains 
a cycle of length at most $2r+3<2r+4$ but the vertex is not within distance
$r+1$ of any such cycle.) Starting with $U=\emptyset$  
let $v_1$ be an arbitrary vertex of $Z$, add it to $U$ and remove 
all vertices of $Z$ within distance  $2r+2$ of $v_1$. 
Clearly at most $d^{2r+2}$ vertices have been deleted. Let
$v_2$ be an arbitrary vertex left in $Z$, add it to $U$ and remove all vertices
of $U$ within distance $2r+2$ of $v_2$. Continuing in this manner
we get a set $U$ of  at least $\frac{n}{2d^{2r+3}}$ vertices.  It is clear
that this set satisfies all the conclusions of the lemma. It is also clear
that $U$ can be computed in polynomial time. \hfill $\Box$

The next lemma about the delocalization of eigenvectors of regular graphs
in cycle-free neighborhood can be proved using the method of Kahale in \cite{Ka}
(see also \cite{AGS} for a recent application of this technique). 
Here we present a simple self contained proof.
\begin{lemma}
\label{l42}
Let $G=(V,E)$ be a $d$-regular graph where $d \geq 3$, 
let $uv$ be an edge of $G$
and suppose that the $r$-neighborhood of $uv$ contains no
cycle. For each $i$ satisfying $0 \leq i \leq r$ let
$N_i$ denote the set of all vertices of distance exactly $i$
from $\{u,v\}$. (In particular, $N_0=\{u,v\}$).
Let $f$ be a nonzero eigenvector of $G$ with eigenvalue
$\mu \geq 2 \sqrt{d-1}$. Then for every $1 \leq i \leq r$
\begin{equation}
\label{e41}
\sum_{w \in N_i} f^2(w) \geq \sum_{w \in N_{i-1}} f^2 (w).
\end{equation}
\end{lemma}
\vspace{0.1cm}

\noindent
{\bf Proof:}\, We apply induction on $i$. Note that by assumption the induced
subgraph of $G$ on the $r$-neighborhood of $uv$ is a $d$-regular tree. Therefore
$|N_i|=2(d-1)^i$ for all $i \leq r$. Let $u_1,u_2, \ldots ,u_{d-1}$
denote the neighbors of $u$ besides $v$, and let $v_1,v_2, \ldots ,v_{d-1}$
denote the neighbors  of $v$ besides $u$. Then 
$$
f(v)+\sum_{i=1}^{d-1} f(u_i)=\mu f(u)
$$
and 
$$
f(u)+\sum_{i=1}^{d-1} f(v_i)=\mu f(v).
$$
By Cauchy-Schwarz,
$$
f^2(v)+\sum_{i=1}^{d-1} f^2(u_i) \geq \frac{\mu^2 f^2(u)}{d} \geq 
\frac{4d-4}{d} f^2 (u)
$$
and similarly
$$
f^2(u)+\sum_{i=1}^{d-1} f^2(v_i) \geq \frac{4d-4}{d} f^2 (v).
$$
Summing, we conclude that
$$
f^2(u)+f^2(v)+\sum_{w \in N_1}f^2(w) \geq \frac{4d-4}{d} (f^2(u)+f^2 (v)),
$$
implying that 
$$
\sum_{w \in N_1}f^2(w) \geq \frac{3d-4}{d} (f^2(u)+f^2 (v))
\geq f^2(u)+f^2(v) =\sum_{w \in N_0} f^2(w).
$$
This proves (\ref{e41}) for $i=1$.

Assuming it holds for $i-1$ we prove it for $i$. For each vertex 
$w \in N_{i-1}$ let $w'$ be its unique parent in $N_{i-2}$ and let
$x_1,x_2, \ldots ,x_{d-1}$ be its neighbors in $N_i$.
Then
$$
f(w')+\sum_{i=1}^{d-1}f(x_i) = \mu f(w).
$$
Since $f(w')=(d-1) \cdot \frac{f(w')}{d-1}$, we get, by Cauchy-Schwarz,
$$
\frac{f^2(w')}{d-1}+\sum_{i=1}^{d-1}f^2(x_i) \geq  \frac{\mu^2 f^2(w)}{2d-2}
\geq 2 f^2(w).
$$
Summing the above inequality  for all $w$ in $N_{i-1}$, each
vertex $w'$ appears in the LHS  exactly $d-1$ times, yielding
$$
\sum_{w' \in N_{i-2}} f^2(w') + \sum_{x \in N_{i}} f^2(x) \geq 
2 \sum_{w \in N_{i-1}} f^2 (w).
$$
This gives
$$
\sum_{x \in N_{i}} f^2(x) \geq  2 \sum_{w \in N_{i-1}} f^2 (w)
- \sum_{w' \in N_{i-2}} f^2(w') \geq \sum_{w \in N_{i-1}} f^2 (w),
$$
where the last inequality follows from the induction hypothesis.
This completes the proof of the induction step, establishing the 
assertion of the lemma. \hfill $\Box$

Finally, we need the main result of Mohanty, O'Donnell and Paredes in
\cite{MOP}, which is the following.
\begin{theo}[\cite{MOP}]
\label{t33}
For every $d$, $\epsilon>0$ and (large) $n$ there is an explicit construction
of an $(n+o(n),d,\lambda)$-graph with $\lambda \leq 2 \sqrt{d-1}+\epsilon$
so that the $s$ neighborhood of any vertex contains at most one cycle,
where $s \geq (\log \log n)^2 $.
\end{theo}

We note that the result is stated in \cite{MOP} without the conclusion 
about the
cycles, but the version above follows from the proof as presented there.

We are now ready to prove Theorem \ref{t13}.
\vspace{0.2cm}

\noindent
{\bf Proof of Theorem \ref{t13}:}\,
Put $r=\lceil 2/\epsilon \rceil$ and $s=2r+4$.
Let $H=(V,E)$ be an explicit 
$(n+u,d,\lambda)$-graph with $u=o(n)$ and $\lambda \leq
2\sqrt{d-1}+\epsilon/2$ in which the $s$-neighborhood of any vertex
contains at most one cycle. Such an $H$ exists by Theorem
\ref{t33}.
By Lemma \ref{l31} we can find efficiently  
a set $U$ of $u$ vertices in $H$ satisfying the
assertion of the lemma. Omit these vertices from the graph 
to get a graph $H'$ 
and add
a matching $M$ between their neighbors retaining the degree of
regularity $d$. Let $G$ denote the resulting graph. Clearly it is
$d$-regular and has $n$ vertices.  Note that the $r$-neighborhood
of any edge $uv$ of the added matching $M$ contains no cycle.
In order to complete the proof it remains to show that every
nontrivial eigenvalue of $G$ has absolute value at most
$2 \sqrt {d-1}+\epsilon$. Let $A_G$ be the adjacency matrix of 
$G$, $A_{H'}$ the adjacency matrix of $H'$ (on the set of vertices
$V$) and $A_M$ the adjacency
matrix of the graph on the set of vertices $V $
whose edges are those
of the matching $M$. Note that $A_G=A_{H'}+A_M$.  Let $\lambda$ be
a nontrivial eigenvalue of $G$ and let $f$ be a corresponding
eigenvector satisfying $\sum_{v \in V} f^2(v)=1$. Then
\begin{equation}
\label{e91}
\lambda=f^t A_G f=f^t A_{H'} f+ f^t A_M f.
\end{equation}
Since $H'$ is an induced subgraph of $H$ and all nontrivial
eigenvalues of $H$ have absolute value at most
$2 \sqrt{d-1}+\epsilon/2$ it follows, by
eigenvalue interlacing, that 
\begin{equation}
\label{e92}
|f^t A_{H'} f| \leq 2 \sqrt{d-1}+\epsilon/2.
\end{equation}
It is also clear that
\begin{equation}
\label{e93}
|f^t A_{M} f| = |2\sum_{uv \in M} f(u)f(v)| \leq \sum_{uv \in M} 
f^2(u)+f^2(v).
\end{equation}
If $|\lambda| \leq 2 \sqrt{d-1}$ there is nothing to prove, we thus
assume that $\lambda \geq 2 \sqrt{d-1}$. 
Since the $r$-neighborhood of any edge of $M$ contains no cycle,
Lemma \ref{l32} implies that for every such edge $uv$,
$$
f^2(u)+f^2(v)  \leq \frac{1}{r} \sum_{w \in N(u,v,r)} f^2(w),
$$
where $N(u,v,r)$ denotes the $r$-neighborhood of $uv$. Since all
these neighborhoods are pairwise disjoint  it follows that
\begin{equation}
\label{e94}
\sum_{uv \in M} f^2(u)+f^2(v) \leq \frac{1}{r} \sum_{w \in V} f^2
(v) \leq \frac{\epsilon}{2}.
\end{equation}
The desired result follows by
plugging (\ref{e92}) and (\ref{e94}) in (\ref{e91}) (using
(\ref{e93})).
\hfill $\Box$

\section{Concluding remarks}

\begin{itemize}
\item
Morgenstern \cite{Mo} gave a strongly explicit construction of
Ramanujan graphs for every degree which is a prime power plus $1$,
but we cannot apply his construction in the proof of Proposition
\ref{p11} since his construction provides Cayley graphs of
different groups (and different sizes) for different degrees and hence 
one cannot pack the graphs corresponding to several degrees.
Similarly, we cannot use his construction in the proof of
Theorem \ref{t12} since for every fixed degree the 
sequence of possible numbers of vertices in his construction 
for this degree is too sparse.
\item
The proof of Theorem \ref{t13} can be applied directly to high girth
Ramanujan graphs like those in \cite{LPS}, \cite{Ma}  in case the 
required degree is $p+1$ for a prime $p$ congruent to 
$1 \bmod 4$ to obtain near Ramanujan graphs of this degree 
with any required
(large) number of vertices.
\item
The problem of obtaining 
strongly explicit (two-sided) Ramanujan (and not nearly Ramanujan) 
graphs for 
any degree and
number of vertices remains open.
\end{itemize}

\noindent
{\bf Acknowledgment}
I thank L\'aszl\'o Babai, Oded Goldreich, Ryan O'Donnell, Ori
Parzanchevski and Peter Sarnak for helpful discussions.

\end{document}